\documentclass[11pt]{article}
\usepackage{graphicx}
\usepackage{amsfonts}
\usepackage{theorem}
\newtheorem{Lem}{Lemma}[section]

\newtheorem{The}[Lem]{Theorem}
\newtheorem{Prop}[Lem]{Proposition}
\newtheorem{Cor}[Lem]{Corollary}

\newtheorem{Rem}[Lem]{Remark}

\newcommand{\qed}{\hbox{\rule{6pt}{6pt}}}

\begin{document}
\title{Inequalities for relative operator entropies and operator means}
\author{Shigeru Furuichi$^1$\footnote{E-mail:furuichi@chs.nihon-u.ac.jp} and Nicu\c{s}or Minculete$^2$\footnote{E-mail:minculeten@yahoo.com}\\
$^1${\small Department of Information Science,}\\
{\small College of Humanities and Sciences, Nihon University,}\\
{\small 3-25-40, Sakurajyousui, Setagaya-ku, Tokyo, 156-8550, Japan}\\
$^2${\small Transilvania University of Bra\c{s}ov, Bra\c{s}ov, 500091, Rom{a}nia}}
\date{}
\maketitle
{\bf Abstract.}
The main purpose of this article is to study estimates for the Tsallis relative operator entropy, by the use of Hermite-Hadamard inequality. Thus, we obtain alternative bounds for the Tsallis relative operator entropy. In the process to derive these bounds, we established the significant relation between the Tsallis relative operator entropy and the generalized relative operator entropy. In addition, we study the properties on monotonicity for the weight of operator means, and for the parameter of relative operator entropies.
\vspace{3mm}


\vspace{3mm}
\vspace{3mm}

\section{Introduction}
In operator theory, we find various characterizations and the relationship 
between operator monotonicity and operator convexity given 
by Hansen and Pedersen \cite{HP}, Chansangiam \cite{Chan}.
In \cite{KA}, Kubo and Ando has studied the connections between 
operator monotone functions and operator means. 
The operator monotone function plays an important roles 
in the theory given by Kubo and Ando. 
Other information about applications of operator monotone functions 
to theory of operator mean can be find in \cite{Ped}. 
Theory of operator mean plays a central role in operator inequalities, 
operator equations, network theory, and quantum information theory. 

Denote by $B(H)$ 
the algebra of bounded linear operators on a  Hilbert space $H$. 
We write $A> 0$ to mean that $A$ is a strictly positive operator, 
or equivalently, $A \geq 0$  and $A$ is invertible.
Furuta and Yanagida showed the following inequality with elegant proof \cite{Furuta}:
\begin{equation}
A{!_p}B \le A{\# _p}B \le A{\nabla _p}B,
\end{equation}
where we respectively denote $p$-weighted harmonic operator mean,  $p$-weighted geometric operator mean and  $p$-weighted arithmetic operator mean by
$A!_pB \equiv \left\{ (1-p)A^{-1} +p B^{-1}\right\}^{-1}$, $A\#_pB \equiv A^{1/2}\left(A^{-1/2}BA^{-1/2}\right)^{p}A^{1/2}$ and $A\nabla_pB \equiv (1-p)A+pB$ for $A, B>0$ and $p \in [0,1]$.
 
On the other hand, Tsallis defined the one-parameter extended entropy for the analysis of a physical model in statistical physics in \cite{Tsa1988}. The properties of the Tsallis relative entropy was studied in \cite{F2006,FYK2004}, by Furuichi, Yanagi and Kuriyama. The relative operator entropy
 \[S\left( {A|B} \right): = {A^{1/2}}\log \left( {{A^{ - 1/2}}B{A^{ - 1/2}}} \right){A^{1/2}}\]
for two invertible positive operators $A$ and $B$ on a Hilbert space, was introduced by Fujii and Kamei in \cite{FK1989}. The parametric extension of the relative operator entropy was introduced by Furuta in \cite{Furuta2004} as
 \[S_p\left( {A|B} \right): = {A^{1/2}} \left(A^{-1/2}BA^{-1/2}\right)^p\log \left( {{A^{ - 1/2}}B{A^{ - 1/2}}} \right){A^{1/2}}\]
for $p \in \mathbb{R}$ and two invertible positive operators $A$ and $B$ on a Hilbert space. Note that $S_0(A|B) \equiv \lim_{p\to 0}S_p(A|B) = S(A|B)$.
In \cite{YKF2005}, Yanagi, Kuriyama and Furuichi introduced a parametric extension of relative operator entropy by the concept of Tsallis relative entropy for operators, as
 \[{T_p}\left( {A|B} \right): = \frac{{{A^{1/2}}{{\left( {{A^{ - 1/2}}B{A^{ - 1/2}}} \right)}^p}{A^{1/2}} - A}}{p},\,\,\,\left( { - 1 \le p \le 1,\,\,p \ne 0} \right)\]
where $A$ and $B$ are two strictly positive operators on a Hilbert space $H$.
In \cite{FYK2005}, we found several results about the Tsallis relative operator entropy. Furuta \cite{Furuta2006} showed two reverse inequalities involving Tsallis relative operator entropy $T_p(A|B)$ via generalized Kantorovich constant $K(p)$.
The Tsallis relative operator entropy can be rewritten as 
 $${T_p}\left( {A|B} \right) = \frac{{A{\natural _p}B - A}}{p},$$
where $A{\natural _p}B: = {A^{1/2}}{\left( {{A^{ - 1/2}}B{A^{ - 1/2}}} \right)^p}{A^{1/2}}$ for all $p \in \mathbb{R}$.
To study Tsallis relative operator entropy is often strongly connected to the study of the $p$-weighted geometric operator mean.
It is known that \cite{FYK2005}:
\begin{equation} \label{ineq_FYK2005}
A - A{B^{ - 1}}A \le {T_p}\left( {A|B} \right) \le B - A,
\end{equation}
for strictly positive operators $A$, $B$ and $p \in \left[ { - 1,0} \right) \cup \left( {0,1} \right]$  and $\mathop {\lim }\limits_{p \to 0} {T_p}\left( {A|B} \right) = S\left( {A|B} \right)$.


\section{Alternative estimate of Tsallis relative operator entropy}
We start from the following known properties of the Tsallis relative operator entropy. See \cite[Theorem 1]{IIKTW2012} or \cite[Theorem 2.5 (ii)]{IIKTW2015} for example. 
\begin{Prop} \label{prop1}
For any strictly positive operators $A$ and $B$ and $p,q \in [-1,0) \cup (0,1]$ with $p \leq q$, we have
$$
T_p(A|B) \leq T_q(A|B).
$$ 
\end{Prop}
This proposition can be proven by the monotone increasing of $\frac{x^p-1}{p}$ on $p \in [-1,0) \cup (0,1]$ for any $x >0$, and implies the following inequalities (which include the inequalities (\ref{ineq_FYK2005})) \cite{YKF2005}:
$$
A-AB^{-1}A=T_{-1}(A|B) \leq T_{-p}(A|B) \leq S(A|B) \leq T_p(A|B) \leq T_1(A|B) =B-A
$$
for any strictly positive operators $A$ and $B$ and $p \in (0,1]$. The general results were recently established in \cite{Nik2014} by the notion of perspective functions. In addition, quite recently the interesting and significant results for relative operator entropy were given in \cite{DB2016} for the case $B \geq A$.
In this section, we treat the relations on the Tsallis relative operator entropy under the assumption such that strictly positive operators $A$ and $B$ have the ordering $A \leq B$ or $A \geq B$. 

In \cite{RFM2017}, we obtained the estimates on Tsallis relative operator entropy by the use of Hermite-Hadamard inequality:
$$
f\left(\frac{a+b}{2} \right) \leq \frac{1}{b-a} \int_a^b f(t) dt \leq \frac{f(a) +f(b)}{2}
$$
for a convex function $f(t)$ defined on the interval $[a,b]$ with $a \neq b$.

\begin{The} {\bf (\cite{RFM2017})} \label{theorem_a}
For any invertible positive operator $A$ and $B$ such that $A\le B$, and $-1 \leq p \le 1$ with $p \neq 0$ we have
\begin{eqnarray*}
&& {{A}^{1/2}}{{\left( \frac{{{A}^{-1/2}}B{{A}^{-1/2}}+I}{2} \right)}^{p-1}}\left( {{A}^{-1/2}}B{{A}^{-1/2}}-I \right){{A}^{1/2}} \\ 
&& \le {{T}_{p}}\left( A|B \right) 
 \le \frac{1}{2}\left( A{{\#}_{p}}B-A{{\natural}_{p-1}}B+B-A \right),  
\end{eqnarray*}
where $I$ is the identity operator.
\end{The}

The inequalities in Theorem \ref{theorem_a} are improvements of the inequalities (\ref{ineq_FYK2005}).
In the present paper, we give the alternative bounds for the Tsallis relative operator entropy. The condition $A \leq B$ in Theorem \ref{theorem_a} can be modified by $uA \le B \le vA$ with $u \ge 1$ so that we use this style (which is often called a sandwich condition) in the present paper. Note that the condition $uA \le B \le vA$ with $u \ge 1$ includes the condition $A \le B$ as a special case, also the condition $uA \le B \le vA$ with $v \le 1$ includes the condition $B \le A$ as a special case.

\begin{The} \label{theorem1}
Let $A$ and $B$ be strictly positive operators such that $uA \le B \le vA$ with $u,v>0$ and let $-1 \le p \le 1$ with $p \neq 0$.
If $u \ge 1$, then
\begin{equation} \label{ineq01_theorem1}
S_{p/2}(A|B) \leq T_p(A|B) \leq \frac{S(A|B) +S_p(A|B)}{2}.
\end{equation}
If $v \le 1$, then the reverse inequalities in (\ref{ineq01_theorem1}) hold.
\end{The}

{\it Proof}:
For $x \geq 1$ and $-1 \leq p \leq 1$ with $p \neq 0$, we define the function $f(t) = x^{pt} \log x$ on $0\leq t \leq 1$. Since $\frac{d^2f(t)}{dt^2} = p^2x^{pt} \left( \log x\right)^3 \geq 0$ for $x \geq 1$, the function $f(t)$ is convex on $t$, for the case $x \geq 1$. Thus we have
\begin{equation} \label{ineq01_proof_theorem1}
x^{p/2} \log x \leq \frac{x^p -1 }{p} \leq \left(\frac{x^p +1 }{2} \right) \log x
\end{equation} 
by Hermite-Hadamard inequality, since $\int_0^1 f(t) dt = \frac{x^p -1}{p}$. 
Note that $I \le uI \le A^{-1/2}BA^{-1/2} \le vI$ from the condition $u \ge 1$. By Kubo-Ando theory \cite{KA}, it is known that for the representing function
$f_m(x)=1mx$ for operator mean $m$, the scalar inequality $f_m(x) \leq f_n(x),\,\,(x>0)$ is equivalent to the operator inequality $AmB \leq AnB$ for all strictly positive operators $A$ and $B$. (Hereafter we omit this description for simplicity in the following proofs.) 
Thus we have the inequality 
\begin{eqnarray*}
&&A^{1/2}\left(A^{-1/2}BA^{-1/2} \right)^{p/2} \log \left(A^{-1/2}BA^{-1/2} \right)A^{1/2}
\leq \frac{A\#_p B-A}{p} \\
&&\leq \frac{A^{1/2}\log \left(A^{-1/2}BA^{-1/2} \right)A^{1/2}+A^{1/2}\left(A^{-1/2}BA^{-1/2} \right)^{p} \log \left(A^{-1/2}BA^{-1/2} \right)A^{1/2}}{2}
\end{eqnarray*}
which is the inequality (\ref{ineq01_theorem1}).
The reverse inequalities for the case $v \leq 1$ can be similarly shown by the concavity of the function $f(t)$ on $t$, for the case $0 < x \leq 1$, taking into account the condition $0 < uI \le A^{-1/2}BA^{-1/2} \le vI \le I$.

\hfill \qed

We note that both side in the inequalities (\ref{ineq01_theorem1}) and their reverses converges to $S(A|B)$ in the limit $p \to 0$. 
From the proof of Theorem \ref{theorem1}, for strictly positive operators $A$ and $B$, we see the following interesting relation between the Tsallis relative operator entropy $T_p(A|B)$ and the generalized relative operator entropy $S_p(A|B)$,
$$
\int_0^1 S_{pt}(A|B) dt =T_p(A|B).
$$

\begin{Rem} \label{remark0}
Let $A$ and $B$ be strictly positive operators such that $uA \le B \le vA$ with $u,v>0$ and let $-1 \le p \le 1$ with $p \neq 0$.
For the case $0< p \leq 1$ and $u \ge 1$, we see
$$
S(A|B) \leq S_{p/2}(A|B) \leq T_p(A|B) \leq \frac{S(A|B) +S_p(A|B)}{2} \leq S_p(A|B)
$$
from the inequalities (\ref{ineq01_theorem1}) since $x^p \log x$ is monotone increasing on $0< p \leq 1$ and $\left( \frac{x^p +1 }{2}\right) \log x \leq x^p \log x$ for $x \geq 1$ and $0 < p \leq 1$.  For the case $-1 \leq p <0$ and $v \le 1$, we also see that the reverse inequalities hold since $x^p \log x$ is monotone increasing on $-1 \leq p < 0$ and $\left( \frac{x^p +1 }{2}\right) \log x \geq x^p \log x$ for $0 < x \leq 1$ and $-1 \leq p < 0$.  
\end{Rem}

\begin{Rem} \label{remark1}
We compare the bounds of  $\frac{{{x^p} - 1}}{p}$ in the inequalities (\ref{ineq01_theorem1}) with the result given in \cite{RFM2017}:
\[{\left( {\frac{{x + 1}}{2}} \right)^{p - 1}}\left( {x - 1} \right) \le \frac{{{x^p} - 1}}{p} \le \left( {\frac{{{x^{p - 1}} + 1}}{2}} \right)\left( {x - 1} \right),\,\,\,\,\left( {x \ge 1,\,\,\,0 < p \le 1} \right).\]
\begin{itemize}
\item[(i)] We have no ordering between $x^{p/2} \log x$ and $\left(\frac{x+1}{2} \right)^{p-1}\left(x-1\right)$. Indeed, when $p=1/4$ and $x=3$, $x^{p/2} \log x -\left(\frac{x+1}{2} \right)^{p-1}\left(x-1\right) \simeq 0.071123.$ On the other hand, when $p=3/4$ and $x=3$, $x^{p/2} \log x -\left(\frac{x+1}{2} \right)^{p-1}\left(x-1\right) \simeq -0.023104.$
\item[(ii)] We have no ordering between $\left(\frac{x^{p}+1}{2}\right) \log x$ and $\left(\frac{x^{p-1}+1}{2} \right)\left(x-1\right)$. Indeed, when $p=1/4$ and $x=3$, $\left(\frac{x^{p-1}+1}{2} \right)\left(x-1\right) - \left(\frac{x^{p}+1}{2}\right) \log x \simeq 0.166458.$ On the other hand, when $p=3/4$ and $x=3$, $\left(\frac{x^{p-1}+1}{2} \right)\left(x-1\right) - \left(\frac{x^{p}+1}{2}\right) \log x \simeq -0.0416177.$
\end{itemize}
Therefore we claim Theorem \ref{theorem1} is not trivial result.
\end{Rem}

\begin{The} \label{theorem2}
Let $A$ and $B$ be strictly positive operators such that $uA \le B \le vA$ with $u \ge 1$ and let $-1 \le p \le 1$ with $p \neq 0$.
Then we have
\begin{eqnarray*}
&&\frac{T_p(A|B)-T_{p-1}(A|B)}{2} \leq 4 \left\{ T_p\left( A\left| \frac{A+B}{2}\right. \right)-T_{p-1}\left( A\left| \frac{A+B}{2}\right. \right)\right\} \\
&&\leq \frac{T_p(A|B)-T_1(A|B)}{p-1} \leq \frac{T_p(A|B)-T_{p-1}(A|B)}{2} +\frac{A \natural_2(B-A)}{4}.
\end{eqnarray*}
\end{The}

{\it Proof}:
It is sufficient to prove the following inequalities for $t \geq 1$ and  $-1 \le p \le 1$ with $p \neq 0$,
\begin{equation} \label{ineq01_proof_theorem2}
l_p(t) \leq k_p(t) \leq c_p(t) \leq l_p(t) +\frac{(t-1)^2}{4},
\end{equation}
where
\begin{eqnarray*}
&& l_p(t) : = \frac{1}{2}\left( \frac{t^p-1}{p} -  \frac{t^{p-1}-1}{p-1}\right),\quad c_p(t) : = \frac{t^p -1}{p(p-1)} -\frac{t-1}{p-1},\\
&& k_p(t) : = \frac{4}{p}\left\{ \left( \frac{t+1}{2}\right)^p -1\right\} -\frac{4}{p-1}\left\{ \left( \frac{t+1}{2}\right)^{p-1} -1\right\}.
\end{eqnarray*}

Firstly, to prove $l_p(t) \leq k_p(t)$, we set the function $h_p(t) : = k_p(t) -l_p(t)$. Then we calculate 
$$
\frac{dh_p(t)}{dt} =(t-1) \left\{ \left( \frac{t+1}{2}\right)^{p-2} -\frac{t^{p-2}}{2}\right\}. 
$$
We set $g_p(t) : = \left( \frac{t+1}{2}\right)^{p-2} -\frac{t^{p-2}}{2}$. Then we have
$$
\frac{dg_p(t)}{dt} =\frac{p-2}{2} \left\{\left( \frac{t+1}{2}\right)^{p-3} -t^{p-3} \right\} \leq 0,\quad g_p(1) = \frac{1}{2}, \quad \lim_{t\to \infty} g_p(t) =0.
$$
Thus we have $g_p(t) \geq 0$, that is, $\frac{dh_p(t)}{dt} \geq 0$ so that we have
$h_p(t) \geq h_p(1)=0$.

Secondly, the inequalities $k_p(t) \leq c_p(t) \leq l_p(t) +\frac{(t-1)^2}{4}$ can be proven in the following way.
We consider $x \geq 1$ and the function $f:[1,x] \to \mathbb{R}$  defined by $f(y) =y^{p-2}$ with $p\in(0,1]$. It follows that $f'\left( y \right) = \left( {p - 2} \right){y^{p - 3}}$ with $f''\left( y \right) = \left( {p - 2} \right)\left( {p - 3} \right){y^{p - 4}} \ge 0$, so the function $f$ is convex. Therefore, we apply the Hermite-Hadamard inequality, we have
\[\left( {x - 1} \right){\left( {\frac{{x + 1}}{2}} \right)^{p - 2}} \le \int\limits_1^x {{y^{p - 2}}dy \le \left( {x - 1} \right)\left(\frac{{{x^{p - 2}} + 1}}{2}\right)} \]
which, by integrating, is equivalent to the inequality
\[\int\limits_1^t {\left( {x - 1} \right){{\left( {\frac{{x + 1}}{2}} \right)}^{p - 2}}dx}  \le \int\limits_1^t {\int\limits_1^x {{y^{p - 2}}dydx} }  \le \int\limits_1^t {\left( {x - 1} \right)\left(\frac{{{x^{p - 2}} + 1}}{2}\right)} dx.\]
Since we have the computations of the following integrals, for $t,x \geq 1$
\[\int\limits_1^t {\int\limits_1^x {{y^{p - 2}}dydx} }  = \frac{{{t^p} - 1}}{{p\left( {p - 1} \right)}} - \frac{{\left( {t - 1} \right)}}{{\left( {p - 1} \right)}},\]
\begin{equation} \label{eq01_theorem1}
\int\limits_1^t {\left( {x - 1} \right){{\left( {\frac{{x + 1}}{2}} \right)}^{p - 2}}dx}  = \frac{2}{{p - 1}}\left( {t - 1} \right){\left( {\frac{{t + 1}}{2}} \right)^{p - 1}} - \frac{4}{{p\left( {p - 1} \right)}}\left\{ {{{\left( {\frac{{t + 1}}{2}} \right)}^p} - 1} \right\}
\end{equation}
and
\begin{equation} \label{eq02_theorem1}
\int\limits_1^t {\left( {x - 1} \right) \left(\frac{{{x^{p - 2}} + 1}}{2}\right)} dx = \frac{1}{2}\left\{ {\frac{{{t^p} - 1}}{p} - \frac{{{t^{p - 1}} - 1}}{{p - 1}} + \frac{1}{2}\left( {{t^2} - 1} \right) - \left( {t - 1} \right)} \right\},
\end{equation}
we obtain the inequality
\begin{eqnarray*} 
&&\frac{2}{{p - 1}}\left( {t - 1} \right){\left( {\frac{{t + 1}}{2}} \right)^{p - 1}} - \frac{4}{{p\left( {p - 1} \right)}}\left\{ {{{\left( {\frac{{t + 1}}{2}} \right)}^p} - 1} \right\}\\
&& \le \frac{{{t^p} - 1}}{{p\left( {p - 1} \right)}} - \frac{{\left( {t - 1} \right)}}{{\left( {p - 1} \right)}} \le \frac{1}{2}\left\{ {\frac{{{t^p} - 1}}{p} - \frac{{{t^{p - 1}} - 1}}{{p - 1}} + \frac{1}{2}\left( {{t^2} - 1} \right) - \left( {t - 1} \right)} \right\}
\end{eqnarray*}
By simple calculations, we find the above inequalities are equivalent to the inequalities $k_p(t) \leq c_p(t) \leq l_p(t) +\frac{(t-1)^2}{4}$.

\hfill \qed

\begin{Rem}
We compare Theorem \ref{theorem2} and Theorem \ref{theorem_a} in \cite{RFM2017}.
The inequalities $ k_p(t) \leq c_p(t) \leq l_p(t) +\frac{(t-1)^2}{4}$ given in (\ref{ineq01_proof_theorem2}) are equivalent to the following inequalities
\begin{equation} \label{ineq01_remark25}
\alpha_p(t) \leq \frac{t^p -1}{p} \leq \beta_p(t)
\end{equation}
where
\begin{eqnarray*} 
&& \alpha_p(t) : = \frac{t^{p-1} -1}{p-3} +\frac{p-1}{2(3-p)}(t^2-1) +(t-1), \\
&& \beta_p(t) : =  (t-1) +2 \left( \frac{t+1}{2}\right)^{p-1} -\frac{4}{p} \left\{ \left( \frac{t+1}{2}\right)^{p} -1\right\}.
\end{eqnarray*}
By the inequalities  (\ref{ineq01_remark25}) with Kubo-Ando theory \cite{KA}, we have
\begin{eqnarray} 
&& B - \frac{3}{2}A - \frac{1 - p}{2\left( 3 - p \right)}BA^{ - 1}B - \frac{1}{3 - p}A\natural_{p-1} B \nonumber \\
&& \le T_p\left(A|B \right)  \le 
B - A + 2\left( BA^{ - 1} - I \right) A \natural_{p-1} \left(\frac{A+B}{2} \right)- 4 T_p \left(A \left| \frac{A+B}{2} \right.\right), \label{ineq02_remark25}
\end{eqnarray} 
which are equivalent to the second and third inequalities given in Theorem \ref{theorem2}.

We compare both bounds of $\frac{t^p -1}{p} $ in (\ref{ineq01_remark25}) with the fundamental inequalities 
\[{\left( {\frac{{t + 1}}{2}} \right)^{p - 1}}\left( {t - 1} \right) \le \frac{{{t^p} - 1}}{p} \le \left( {\frac{{{t^{p - 1}} + 1}}{2}} \right)\left( {t - 1} \right),\,\,\,\,\left( {t \ge 1,\,\,\,0 < p \le 1} \right)\]
to obtain Theorem \ref{theorem_a}. For this purpose, let $t \geq 1$ and $-1 \leq p \leq 1$ with $p\neq 0$. And we set the functions $\gamma_p(t)$ and $\delta_p(t)$ by
$$
\gamma_p(t) : =  \alpha_p(t)  -   {\left( {\frac{{t + 1}}{2}} \right)^{p - 1}}\left( {t - 1} \right), \quad \delta_p(t) : = \left( {\frac{{{t^{p - 1}} + 1}}{2}} \right)\left( {t - 1} \right) - \beta_p(t).
$$
By numerical computations, we have the following results.
\begin{itemize}
\item[(i)] $\gamma_{1/2}(3/2) \simeq 0.00118777$ and $\gamma_{1/2}(5/2) \simeq -0.0118756$.
\item[(ii)] $\delta_{1/2}(3/2) \simeq -0.890458$ and $\delta_{1/2}(5/2) \simeq 0.795489$.
\end{itemize}
Thus we conclude that for the Tsallis relative entropy $T_p(A|B)$ there is no ordering between the bounds given in  the inequalities  (\ref{ineq02_remark25})  and ones given in Theorem \ref{theorem_a} of \cite{RFM2017}.
Therefore we claim Theorem \ref{theorem2} is also not trivial result.
\end{Rem}

Taking the limit $p \to 0$ in Theorem \ref{theorem2}, we have the following corollary.
\begin{Cor}
For strictly positive operators $A$ and $B$ such that $uA \le B \le vA$ with $u \ge 1$, we have
\begin{eqnarray*}
&&\frac{S(A|B)-T_{-1}(A|B)}{2} \leq 4 \left\{ S\left( A\left| \frac{A+B}{2}\right. \right)-T_{-1}\left( A\left| \frac{A+B}{2}\right. \right)\right\} \\
&& \leq T_1(A|B) - S(A|B) \leq \frac{S(A|B)-T_{-1}(A|B)}{2} +\frac{A \natural_2(B-A)}{4}.
\end{eqnarray*}
\end{Cor}


\section{Monotonicity on the parameter of relative operator entropies}

In our previous section, we gave the interesting relations between Tsallis relative operator entropy $T_p(A|B)$ \cite{YKF2005} and the generalized relative operator entropy $S_p(A|B)$ \cite{Furuta2004}.
In this section, we study the monotonicity on parameter $p$ related to two relative operator entropies $T_p(A|B)$ and $S_p(A|B)$.

\begin{Lem} \label{lemma01_section3}
For $x \geq 1$, we have the inequality $1-x+x \log x \leq x\left( \log x\right)^2$.
For $\frac{1}{e} \leq x \leq 1$, we also have the same inequality.
\end{Lem} 

{\it Proof}:
For $x \geq 1$, we set the function $g(x)\equiv 1-x+x\log x-x \left( \log x\right)^2$. Then we calculate $g'(x)=-(1+ \log x)\log x \leq 0$ so that $g(x) \leq g(1) =0$ for $x \geq 1$. We also have $g'(x) \geq 0$ for $\frac{1}{e}\leq x \leq 1$ so that we have $g(x) \leq g(1) =0$.

\hfill \qed

\begin{Prop}\label{prop_3_2}
Let $A$ and $B$ be strictly positive operators such that $uA \le B \le vA$ with $u,v>0$ and let $-1 \le p \le 1$ with $p \neq 0$.
If we have the condition either (i) $u \ge 1$ and $0< p \leq q \leq 1$ or
(ii) $v \leq 1$ and $-1 \leq p \leq q <0$, then
$$
T_p(A|B)-S_p(A|B) \geq T_q(A|B)-S_q(A|B). 
$$
If we also have the condition either (iii) $e^{-1/q} \leq v \leq 1$ and $0< p \leq q \leq 1$ or
(iv) $1 \leq u \leq e^{-1/p}$ and $-1 \leq p \leq q <0$, then the above inequality holds.
\end{Prop}

{\it Proof}:
For $t>0$ and $-1\leq p \leq 1$ with $p \neq 0$, we set the function
$f(p,t) \equiv \frac{t^p-1}{p} -t^p\log t$. Then we  calculate
$\frac{df(p,t)}{dp} =\frac{1}{p^2}\left(1-t^p+t^p\log t^p- t^p \left( \log t^p\right)^2\right) \leq 0$. By the use of Lemma \ref{lemma01_section3}
with $x \equiv t^p \geq 1$ for both cases (i)$t\geq 1$ and $0<p\leq 1$ or (ii)$t \leq 1$ and $-1 \leq p <0$, the desired inequality holds.
From Lemma \ref{lemma01_section3}, we also find that $\frac{df(p,t)}{dp}\leq 0$ for $\frac{1}{e}\leq t^p \leq 1$ so that the desired inequality holds for both cases  (iii) $e^{-1/q} \leq v \leq 1$ and $0< p \leq q \leq 1$ or
(iv) $1 \leq u \leq e^{-1/p}$ and $-1 \leq p \leq q <0$.

\hfill \qed

\begin{Lem}\label{lemma02_section3}
Define $g(x) \equiv 1-x+x\log x -\frac{1}{2}x (\log x)^2$ for $x >0$.
If $0 < x \leq 1$, then $g(x) \geq 0$. If $x \geq 1$, then $g(x) \leq 0$.
\end{Lem}

{\it Proof}:
It is trivial from $\frac{dg(x)}{dx} = -\frac{1}{2}(\log x)^2$.

\hfill \qed

\begin{Prop}\label{prop_3_4}
Let $A$ and $B$ be strictly positive operators such that $uA \le B \le vA$ with $u,v>0$ and let $-1 \le p \le 1$ with $p \neq 0$.
If we have the condition either (i) $u \ge 1$ and $-1 \leq p \leq q <0$ or
(ii) $v \leq 1$ and $0<  p \leq q \leq 1$, then
$$
T_p(A|B)-\frac{1}{2}S_p(A|B) \leq T_q(A|B)-\frac{1}{2}S_q(A|B). 
$$
If we have the condition either (iii) $u \geq 1$ and $0<  p \leq q \leq 1$ or
(iv) $v \leq 1$ and $-1 \leq p \leq q <0$, then
$$
T_p(A|B)-\frac{1}{2}S_p(A|B) \geq T_q(A|B)-\frac{1}{2}S_q(A|B). 
$$
\end{Prop}

{\it Proof}:
We set the function $f(p,t) \equiv \frac{t^p-1}{p}-\frac{1}{2}t^p\log t$ for $t >0$ and $-1 \leq p \leq 1$ with $p \neq 0$.
Then we calculate 
$\frac{df(p,t)}{dp} = \frac{1}{p^2}\left( 1-t^p+t^p\log t^p-\frac{1}{2}t^p\left( \log t^p\right)^2\right)$.
From Lemma \ref{lemma02_section3} with $x \equiv t^p$, we find $\frac{df(p,t)}{dp} \geq 0$ under the condition either (i) $t \geq 1$ and $-1\leq p <0$ or (ii) $0<t \leq 1$ and $0<p \leq 1$. Similarly from  Lemma \ref{lemma02_section3} with $x \equiv t^p$, we find $\frac{df(p,t)}{dp} \leq 0$ under the condition either (iii) $t \geq 1$ and $0< p \leq 1$ or (iv) $0<t \leq 1$ and $-1 \leq p <0$. These imply the conclusion of this proposition, by Kubo-Ando theory.

\hfill \qed

Comparing Proposition \ref{prop_3_2} and Proposition \ref{prop_3_4}, we show slightly precise results, by the similar way to these propositions.
For this purpose, we prepare the following lemma.

\begin{Lem}\label{lemma_3_5}
For $x>0$ and $c \in \mathbb{R}$, we set the function $g(x)\equiv 1-x+x\log x-cx(\log x)^2$.
Then we have $g(x) \geq 0$ under the following three conditions (a) $0<x \leq 1$ and $0<c\leq \frac{1}{2}$, (b)
$1 \leq x \leq e^{\frac{1-2c}{c}}$ and $0<c\leq \frac{1}{2}$, or (c) $x >0$ and $c \leq 0$. We also have
$g(x) \leq 0$ under the following two conditions (d) $e^{\frac{1-2c}{c}} \leq x \leq 1$ and $\frac{1}{2} \leq c$ or (e) $x \geq 1$ and $\frac{1}{2} \leq c$.
\end{Lem}

{\it Proof}:
Since we have $1-x+x\log x \geq 0$ for $x >0$, we have $g(x) \geq 0$ for the case (c). From here we assume $c \neq 0$.
We calculate $g'(x)=(1-2c-c\log x) \log x$. Then we easily have $g(x) \geq 0$ for the case (a), and $g(x) \leq 0$ for the case (e). As for the case (b), we find $g'(x) \geq 0$ for $1 \leq x \leq e^{\frac{1-2c}{c}}$
so that $g(x) \geq g(1)=0$ for the case (b). As for the case (d), we also find $g'(x) \geq 0$ for $ e^{\frac{1-2c}{c}} \leq x \leq 1$
so that $g(x) \leq g(1)=0$ for the case (d). 

\hfill \qed

Note that $l(c) \equiv g\left(e^{\frac{1-2c}{c}} \right) =1+(1-4c)e^{\frac{1-2c}{c}}$ for $c \neq 0$, then
$l'(c)=-\left(\frac{1-2c}{c}\right)^2e^{\frac{1-2c}{c}} \leq 0$. Thus we have $l(c)\geq l\left( \frac{1}{2}\right) =0$ for $c \leq \frac{1}{2}$, and $l(c)\leq l\left( \frac{1}{2}\right) =0$ for $\frac{1}{2}\leq c$.

\begin{Prop}\label{prop_3_6}
Let $A$ and $B$ be strictly positive operators such that $uA \le B \le vA$ with $u,v>0$, $c \in \mathbb{R}$ and let $-1 \le p \le 1$ with $p \neq 0$.
\begin{itemize}
\item[(A)] For $0<c \leq \frac{1}{2}$, we have the inequality
\begin{equation}\label{ineq01_prop_3_6}
T_p(A|B)-c S_p(A|B) \leq T_q(A|B)-c S_q(A|B),
\end{equation}
under the following conditions (a1), (a2), (b1) or (b2).
\begin{itemize}
\item[(a1)] $u \geq 1$ and $-1 \leq p \leq q <0$.
\item[(a2)] $v \leq 1$ and $0 < p \leq q \leq 1$.
\item[(b1)] $1 \leq u \leq v \leq e^{\frac{1-2c}{c q}} \left( \leq e^{\frac{1-2c}{c p}}\right)$ and $0 < p \leq q \leq 1$.
\item[(b2)] $\left(e^{\frac{1-2c}{c q}} \leq \right) e^{\frac{1-2c}{c p}}\leq u  \leq v \leq 1$  and $-1 \leq p \leq q <0$.
\end{itemize}

\item[(B)] For $c \geq \frac{1}{2}$, we have the inequality
$$
T_p(A|B)-c S_p(A|B) \geq T_q(A|B)-c S_q(A|B),
$$
under the following conditions (d1), (d2), (e1) or (e2).
\begin{itemize}
\item[(d1)] $\left(e^{\frac{1-2c}{c p}} \leq \right) e^{\frac{1-2c}{c q}} \leq u \leq v \leq 1$ and $0 < p \leq q \leq 1$.
\item[(d2)] $1 \leq u \leq v \leq e^{\frac{1-2c}{c p}} \left( \leq  e^{\frac{1-2c}{c q}} \right)$ and $-1 \leq p \leq q <0$.
\item[(e1)] $u \geq 1$ and $0 < p \leq q \leq 1$.
\item[(e2)] $v \leq 1$  and $-1 \leq p \leq q <0$.
\end{itemize}
\item[(C)]For $c \leq 0$ and $-1 \leq p \leq q \leq 1$ with $p \neq 0$, $q \neq 0$, we have the inequality (\ref{ineq01_prop_3_6}).
 
\end{itemize}
\end{Prop}

{\it Proof}:
We set $f(p,t)\equiv \frac{t^p-1}{p} -ct^p\log t$ for $t>0$ and $-1 \leq p \leq 1$ with $p \neq 0$.
We calculate $\frac{df(p,t)}{dp} =\frac{1}{p^2}\left(1-t^p+t^p\log t^p-c t^p\left(\log t^p\right)^2\right)$.
From (a), (b) in Lemma \ref{lemma_3_5}, we have $\frac{df(p,t)}{dp}  \geq 0$ under the conditions
(a1),(a2),(b1),(b2) or (c). From (d), (e)  in Lemma \ref{lemma_3_5}, we also have $\frac{df(p,t)}{dp}  \leq 0$ under the conditions (d1), (d2), (e1) or (e2). Finally, from (c)  in Lemma \ref{lemma_3_5}, we have $\frac{df(p,t)}{dp}  \geq 0$ under the conditions (C).
Therefore we have the inequalities in the present proposition.

\hfill \qed


\section{Monotonicity on the weight of operator means}

In this section, along to the previous section, we study the monotonicity of the weight $p$ in weighted mean, since geometric operator mean is used in the definition of Tsallis relative operator entropy.
We review that the following inequalities showing the ordering among three $p$-weighted  means. 
$$
\left\{ (1-p) + pt^{-1} \right\}^{-1} \leq t^p \leq (1-p) + p t,\quad (t>0,\,\, 0\leq p \leq 1).
$$
 We here give the following propositions.

\begin{Prop}\label{proposition2_1}
Let $A$ and $B$ be strictly positive operators, and let $p,q \in (0,1]$.
If $p \leq q$, then
$$\frac{A\nabla_p B -A\sharp_p B}{p} \geq \frac{A\nabla_q B -A\sharp_q B}{q}.$$
\end{Prop}

{\it Proof}:
Since $h(p) : = \frac{x^p - 1}{p} -(x-1)$ is increasing function of $p$ for any $x>0$, if $p \leq q$, then $h(p) \leq h(q)$ which is
\begin{equation} \label{ineq01_prop01}
\frac{x^p -1 -p(x-1)}{p} \leq \frac{x^q -1 -q(x-1)}{q}. 
\end{equation}
By Kubo-Ando theory \cite{KA}, we thus have the desired result.

\hfill \qed

We can obtain the following results in relation to Proposition \ref{proposition2_1}.

\begin{The}
Let $A$ and $B$ be strictly positive operators such that $uA \le B \le vA$ with $u,v>0$ and let $p,q \in (0,1)$ with  $p \leq q$.
If $v \le 1$, then
\begin{equation} \label{ineq01_theorem4_2}
\frac{A\nabla_p B -A\sharp_p B}{p(1-p)} \leq \frac{A\nabla_q B -A\sharp_q B}{q(1-q)}.
\end{equation}
If  $u \ge 1$, then the reverse inequality in (\ref{ineq01_theorem4_2}) holds.
\end{The}

{\it Proof}:
Since $y^{p-2} \geq y^{q-2}$ for $y \leq 1$ and $0<p\leq q <1$, we have
$$
0 \leq \int\limits_x^1 {\int\limits_t^1 {\left( {{y^{p - 2}} - {y^{q - 2}}} \right)dydt = } } \frac{{\left( {1 - p} \right) + px - {x^p}}}{{p\left( {1 - p} \right)}} - \frac{{\left( {1 - q} \right) + qx - {x^q}}}{{q\left( {1 - q} \right)}}.
$$
Thus we have the desired result by Kubo-Ando theory.
Since $y^{p-2} \leq y^{q-2}$ for $y \geq 1$ and $0<p\leq q <1$, in addition $\int\limits_1^x {\int\limits_1^t {{y^{p - 2}}dydt = } } \int\limits_x^1 {\int\limits_t^1 {{y^{p - 2}}dydt} } $ we similarly obtain the following opposite inequality
$$
0 \geq \int\limits_1^x {\int\limits_1^t {\left( {{y^{p - 2}} - {y^{q - 2}}} \right)dydt = } } \frac{{\left( {1 - p} \right) + px - {x^p}}}{{p\left( {1 - p} \right)}} - \frac{{\left( {1 - q} \right) + qx - {x^q}}}{{q\left( {1 - q} \right)}}.
$$
which implies the desired result.

\hfill \qed

\begin{Prop}
Let $A$ and $B$ be strictly positive operators such that $uA \le B \leq vA$ with $v \leq 1$ and let $p,q \in (0,1]$.
If $p \leq q$, then
$$
\frac{A\sharp_p B-A!_p B}{p} \geq \frac{A\sharp_q B-A!_q B}{q} 
$$
\end{Prop}

{\it Proof}:
For $0 < p \leq 1$ and $0 < x \leq 1$, we set $f(p,x) = f_1(p,x)f_2(p,x)$ with
$f_1(p,x) \equiv \frac{x^p}{(1-p) +px^{-1}}$ and $f_2(p,x) \equiv \frac{1}{p}\left( 1-p+px^{-1}-x^{-p}\right)$. Since $(1-p) +pt -t^p \geq 0$ for $t>0$ and $0 \leq p \leq 1$, $f_2(p,x) \geq 0$ for $0 < p \leq 1$ and $0 < x \leq 1$. Putting $t=\frac{1}{x}$ in the inequality (\ref{ineq01_prop01}), we find that $f_2(p,x)$ is decreasing. Since it is trivial that $f_1(p,x) \geq 0$, we finally show $\frac{df_1(p,x)}{dp} \leq 0$. We calculate the first derivative of the function $f_1(p,x)$ by $p$ as
$$
\frac{df_1(p,x)}{dp} = \frac{x^p \left\{ \left(1+\left(x^{-1} -1\right)p\right)\log x -\left( x^{-1} -1\right) \right\} }{ \left\{ (1-p)+px^{-1} \right\}^2}.
$$
Putting $s \equiv x^{-1} -1 \geq 0$, we have
$$
\left(1+\left(x^{-1} -1\right)p\right)\log x -\left( x^{-1} -1\right) =-(ps+1)\log(s+1) -s \leq 0
$$
which implies $\frac{df_1(p,x)}{dp} \leq 0$. Thus we find that  $\frac{df(p,x)}{dp} \leq 0$ for  $0 < p \leq 1$ and $0 < x \leq 1$. Therefore, if $p \leq q$, then
$f(p,x) \geq f(q,x)$.

\hfill \qed

%
%

\begin{Lem} \label{lemma01}
Let $t >0$ and $0 \leq p \leq 1$.
If we have the condition either {\rm (i)} $0 < t \leq 1$ and $0 \leq p \leq \frac{1}{2}$ or {\rm (ii)} $t \geq 1$ and $\frac{1}{2} \leq p \leq 1$, then $(1-p)+pt \geq \frac{t-1}{\log t}$.
\end{Lem}

{\it Proof}:
Both cases are easily proven from 
$$
(1-p) + pt \geq \frac{t+1}{2} \geq \frac{t-1}{\log t}.
$$
The first inequality is true by $\left(p-\frac{1}{2}\right)(t-1) \geq 0$ and the second inequality holds for $t>0$.

\hfill \qed

\begin{Lem}  \label{lemma02}
Let $x >0$ and $0 \leq p \leq 1$.
If  $x \geq 1$ and $0 \leq p \leq \frac{1}{2}$, then
$\log x \geq \frac{x-1}{(1-p)x+p} \geq 0$. If $0 <  x \leq 1$ and $\frac{1}{2} \leq t \leq 1$, then $\log x \leq \frac{x-1}{(1-p)x+p} \leq 0$.
\end{Lem}

{\it Proof}:
Put $x = \frac{1}{t}$ in Lemma \ref{lemma01}.

\hfill \qed

\begin{The}
Let $A$ and $B$ be strictly positive operators such that $uA \le B \le vA$ with $u,v > 0$ and let $p,q \in (0,1]$.
If we have the condition either {\rm (i)} $u \ge 1$ and $0 \leq p \leq q \leq \frac{1}{2}$ or {\rm (ii)} $v \le 1$ and $\frac{1}{2} \leq p \leq q \leq 1$, then 
$$
\frac{A \sharp_p B -A!_pB}{p} +pA^{1/2}\left(\log A^{-1/2}BA^{-1/2}\right)^2A^{1/2} \leq \frac{A \sharp_q B -A!_qB}{q} +qA^{1/2}\left(\log A^{-1/2}BA^{-1/2}\right)^2A^{1/2}. 
$$
\end{The}

{\it Proof}:
We consider the function
$$
f(x,p) \equiv \frac{1}{p}\left(x^p-\frac{x}{(1-p)x+p} \right) +p (\log x)^2.
$$
Then we calculate
\begin{eqnarray*}
\frac{df(x,p)}{dp} &=&\frac{d}{dp}\left( \frac{x^p-1}{p} -\frac{x-1}{(1-p)x +p}\right) + (\log x)^2\\
&=&\frac{d}{dp}\left( \frac{x^p-1}{p}\right) +(\log x)^2 -\left( \frac{x-1}{(1-p)x+p}\right)^2 \geq 0
\end{eqnarray*}
The last inequality is due to Lemma \ref{lemma02} and the fact $\frac{d}{dp}\left(\frac{x^p-1}{p}\right) = \frac{x^p}{p^2} \left(\log x^p -1 +\frac{1}{x^p} \right)\geq 0$ by $\log t \leq t-1$ for $t >0$. Thus we have
$f(x,p) \leq f(x,q)$ under the condition
either {\rm (i)} $x \geq 1$ and $0 \leq p \leq q \leq \frac{1}{2}$ or {\rm (ii)} $0<x \leq 1$ and $\frac{1}{2} \leq p \leq q \leq 1$.

\hfill \qed


\section*{Acknowledgement}
The authors thank anonymous referees for giving valuable comments and suggestions to improve our manuscript. 
The author (S.F.) was partially supported by JSPS KAKENHI Grant Number
16K05257.

\end{document}